\documentclass[11pt,reqno]{amsart}

\usepackage{amssymb}

\usepackage[T1]{fontenc}
\usepackage[dvipsnames]{xcolor}
\usepackage{palatino}
\input amssym.def
\usepackage{amsmath, amsfonts}
\usepackage{amscd, mathtools}
\usepackage[mathscr]{eucal}
\usepackage{palatino}
\newfont{\cyrr}{wncyr10}

\newcommand{\N}{{\mathfrak N}}
\newcommand{\Z}{{\mathbb Z}}
\newcommand{\Q}{{\mathbb Q}}

\newcommand{\rO}{{\rm O}}
\newcommand{\ro}{{\rm o}}

\newcommand{\K}{{\mathbf K}}

\newcommand{\thmref}[1]{Theorem~\ref{#1}}
\newtheorem{thm}{Theorem}[section]

\newtheorem{lem}[thm]{Lemma}
\newtheorem{cor}[thm]{Corollary}

\newcommand{\corref}[1]{Corollary~\ref{#1}}

\newcommand{\lemref}[1]{Lemma~\ref{#1}}

\renewcommand{\b}{{\mathfrak{b}}}
\renewcommand{\a}{{\mathfrak{a}}}

\renewcommand{\P}{{\mathfrak{p}}}

\renewcommand{\O}{{{\mathcal{O}}}}

\begin{document}

\title[Generalized Euler-Kronecker constants]{Transcendence of  
generalized Euler-Kronecker constants}

\author{Neelam Kandhil and Rashi Lunia}

\address{ Neelam Kandhil \\ \newline
	Max-Planck-Institut f\"ur Mathematik,
	Vivatsgasse 7, D-53111 Bonn, Germany.}
\email{kandhil@mpim-bonn.mpg.de}

\address{ Rashi Lunia \\ \newline
The Institute of Mathematical Sciences, 
A CI of Homi Bhabha National Institute, 
CIT Campus, Taramani, Chennai 600 113, India.}
\email{rashisl@imsc.res.in}

\subjclass[2020]{11J81, 11J86}
	
\keywords{Linear forms in logarithms, generalized Euler-Kronecker 
constants.}

\begin{abstract} 
We introduce some generalizations of the Euler-Kronecker constant of a 
number field and study their arithmetic nature.
\end{abstract} 
\maketitle   
\section{\large{Introduction and Preliminaries}}
In 1740, Euler introduced the Euler-Mascheroni constant, which is defined as 
\begin{equation}\label{gamma}
\gamma  = \lim_{x \to \infty} \left( \sum_{n \le x } \frac{1}{n} - \log x \right).
\end{equation}
Since then it has been extensively studied but many questions about its behaviour 
are yet unanswered. 
For example, it is not known if $\gamma$ is rational
or irrational. For a survey on Euler's constant, 
please refer \cite{lag}. In \cite{DF}, Diamond and Ford introduced
a generalization of Euler's constant as follows.
For a non-empty finite set of distinct primes $\Omega$, let 
$P_{\Omega}$ denote the product of the elements of $\Omega$ and $\delta_{\Omega}=
\displaystyle \prod_{\mathfrak{p} \in \Omega} ( 1- \frac{1}{p})$.
Then the generalized Euler constant is defined as
$$
\gamma(\Omega)  = 
\lim_{ x \to \infty} \left(\sum_{ \substack{\atop{n \leq x\atop (n, P_{\Omega})=1}}}
\frac{1}{n}  ~-~  \delta_{\Omega} \log x \right).
$$
Note that, when $\Omega = \emptyset$, we have 
$P_{\Omega}=1=\delta_{\Omega}$ and $\gamma(\Omega) = \gamma $. 
In this context, Murty and Zaytseva \cite{MZ} proved the following theorem.

\begin{thm} {\rm(Murty and Zaytseva)}
At most one number in the infinite list 
$ \{\gamma(\Omega)\}$, as $\Omega$ varies over all finite subsets of 
distinct primes, is algebraic.
\end{thm}

We note that $\gamma$ appears as the
constant term in the Laurent series expansion of $\zeta(s)$ around $s=1$.
This observation led Ihara \cite{YI2} to define the Euler-Kronecker constant
associated to a number field as follows. 

\noindent Let $\K$ be a number field of degree $n$ and let $\O_\K$ denote its ring 
of integers. The Dedekind zeta function of $\K$ is given by 
$$
\zeta_{\K}(s) = \sum_{(0)\neq \a \subseteq{\O_{\K}}} \frac{1}{\N\a^s}, 
\phantom{mm} \Re(s) > 1.
$$
It has a meromorphic continuation to the entire complex plane
with only a simple pole at the point $s=1$. Its Laurent series
expansion around $s=1$ is given by
\begin{eqnarray}\label{Laurentseries}
\zeta_\K(s) =
\frac{\rho_\K}{s -1} + c_\K + \rO(s-1),
\end{eqnarray}
where $\rho_\K \ne 0$ is the residue of $\zeta_\K$ at $s=1.$
The ratio 
\begin{equation}\label{EKdef1}
\gamma_\K := c_\K/\rho_\K,
\end{equation}
was defined by Ihara \cite{YI2} as the Euler-Kronecker constant of $\K$. 
In the next section, an expression analogous to \eqref{gamma}
is given for $\gamma_\K$.

The aim of this article is to study the arithmetic nature of generalizations 
of Euler-Kronecker constants. To do so, we introduce some notations.
Let $\mathcal{P}_{\K}$ denote the set of non-zero prime ideals $\mathfrak{p}$ of $\O_{\K}$ 
and let $\Omega$ be a non-empty subset of $\mathcal{P}_\K$ (possibly infinite) such that
\begin{equation}\label{conv}
\sum_{\mathfrak{p}\in \Omega} \frac {\log \N(\mathfrak{p})}{\N(\mathfrak{p})- 1}  
< \infty.
\end{equation}
For $\K = \mathbb{Q}$, the set of Pjateckii-$\check{S}$apiro primes is an example
of such an infinite subset.
Let  $N_{\Omega}=\{ \mathfrak{p}\cap \Z \mid \mathfrak{p} \in \Omega\}$. We set
\begin{eqnarray}
P(\Omega(x)) 
= \prod_{\mathfrak{p} \in \Omega(x)} \mathfrak{p} 
\phantom{mm}
\text{ and }
\phantom{mm}
\delta_\K(\Omega(x)) 
= \prod_{\mathfrak{p} \in \Omega(x)} \left( 1- \frac{1}{\N(\mathfrak{p})}\right),
\end{eqnarray}
where $\Omega(x) =\{\mathfrak{p} \in \Omega ~| ~ \N(\mathfrak{p})\leq x\}.$
Then by \eqref{conv}, 
$\displaystyle \lim_{x \to \infty}\delta_\K(\Omega(x))$ exists and equals
\begin{equation}
\delta_\K(\Omega)=\prod_{\mathfrak{p}\in \Omega}\left(1-\frac{1}{\N(\mathfrak{p})}\right).
\end{equation}
Note that $\delta_\K(\Omega)=1$ for $\Omega = \emptyset$.
The generalized Euler-Kronecker constant associated to $\Omega$ is denoted by 
$\gamma_\K(\Omega)$ and is defined as
$$
\lim_{ x \to \infty} \left( \frac{1}{\rho_\K} \sum_{ \substack{0\ne \a \subset \O_{\K}
\atop{ \N(\a) \leq x\atop  (\a, P(\Omega(x))=1}}}
\frac{1}{\N(I)}  ~-~  \delta_\K(\Omega(x)) \log x\right).
$$
In \S 3, we will show that this limit exists.
We note that $\gamma_\K(\Omega) = \gamma_\K$ when $\Omega=\emptyset$. In this set up, we have the following
theorem.
\begin{thm}\label{trans}
Let $\{\Omega_i\}_{i\in I}$ be a family of subsets of $\mathcal{P}_\K$ satisfying equation 
\eqref{conv}. Further suppose that $
N_{\Omega_i} \backslash N_{\Omega_j}$ is nonempty and finite for all $ i,j \in I$
and $i \ne j$ . Then at most one number from the infinite list 
$$
\left\{ \frac{\gamma_\K(\Omega_i)}{ \delta_\K(\Omega_i) } ~  \Big | ~ i \in I \right\}
$$
is algebraic. 
\end{thm}

\noindent We digress here a little to make an interesting observation. For $K=\Q$, it 
is known by Merten's theorem that as $x \to \infty$,
$$
\delta_\Q(\Omega(x))\sim \frac{e^{-\gamma}}{\log{x}}.
$$
This makes one wonder if $\gamma_\K$ appears as an exponent in the expression for  
$\K\neq \Q$ too. A result of Rosen \cite{MiRo} shows that this is not true in general.
 More precisely he showed that as $x \to \infty$,
$$
\delta_\K(\Omega(x))\sim \frac{e^{-\gamma}}{\rho_\K \log{x}}.
$$

\section{\large{Preliminaries and Lemmas}}
Let $\K$ be a number field of degree $n$. Throughout this section, $\P$ denotes a 
non-zero prime ideal of $\O_\K$. We recall the following result on counting the number 
of integral ideals of $\O_{\K}$.  

\begin{lem}\cite[Ch 11]{ME}\label{air}
Let $a_m$ be the number of integral ideals of $\mathcal{O}_\K$ with norm $m$. Then, as 
$x$ tends to infinity, 
$$
\sum_{m=1}^{x} a_m = \rho_{\K} x + \rO(x^{1-1/n}).
$$
\end{lem}

Using this result we get the following expression for $\gamma_\K$, which is analogous 
to \eqref{gamma}.

\begin{lem}\label{EKdefn}
For any number field $\K$, the limit
\begin{equation}\label{Abel}
\lim_{x \rightarrow \infty} \left(\frac{1}{\rho_\K}\sum_{0\ne \a \subset \O_{\K} 
\atop \N(\a) \leq x}  \frac{1}{\N(\a)} -  \log x   \right)
\end{equation}
exists and equals $\gamma_{\K}$.
\end{lem}

\begin{proof}
Applying partial summation formula and \lemref{air}, the result follows.
\end{proof}

The M\"obius function $\mu_\K$ and the von Mangoldt function 
$\Lambda_\K$ are defined on $\O_\K$ as follows.

\begin{eqnarray*}
\mu_\K(\a)
=
\begin{cases}
1
& \mbox{ if   }~ \a=\O_\K\\
(-1)^r  & \mbox{ if } \a  \mbox{ is a product of } r 
 \mbox{ distinct prime ideals } \\
0 & \mbox{ otherwise},
\end{cases} 
\end{eqnarray*}

\begin{eqnarray*}
\Lambda_\K(\a) 
= \begin{cases}\log\N\P & \mbox{ if}~ \a=\P^m \mbox{ for some } \P \mbox{ and 
some integer } m \geq 1\\
0  & \mbox{ otherwise}.
\end{cases}  
\end{eqnarray*}

We record the following identities satisfied by these functions which can be derived
using techniques similar to \cite[Exercise: 1.1.2, 1.1.4, 1.1.6]{RM}.
\begin{equation}
\sum_{J \mid I} \frac{\mu_\K(J)}{\N(J)} = \prod_{\mathfrak{p}\mid I} 
\left( 1- \frac{1}{\N(\mathfrak{p})}\right),
\end{equation}

\begin{equation}
\mu_\K(I) \log \N(I) = -\sum_{J\mid I} \Lambda_\K(J) \mu_\K\left(IJ^{-1}\right).
\end{equation}

We end this section by stating the key ingredient in the proof of \thmref{trans}.
\begin{lem}[Lindemann \cite{FL}]\label{FL}
If $\alpha \neq 0, 1$ is an algebraic number, then $\log{\alpha}$ is transcendental,
where $\log$ denotes any branch of logarithmic function. 
\end{lem}

\section{Generalized Euler-Kronecker constants}

Let $\mathcal{P}_\K$ denote the set of non-zero prime ideals of $\O_{\K}$.
For any non-empty finite set $\Omega_f \subset \mathcal{P}_\K$, we set
\begin{eqnarray}
P(\Omega_f)= \prod_{\mathfrak{p} \in \Omega_f} \mathfrak{p} 
\text{ \phantom{m} and \phantom{m}}
\delta_\K(\Omega_f)=\prod_{\mathfrak{p}\in \Omega_f}\left(1-\frac{1}{\N(\mathfrak{p})}
\right),
\end{eqnarray}
with the convention that $P(\Omega_f)=1=\delta_\K(\Omega_f)$, when $\Omega_f=\emptyset$. Since $\mathcal{O}_\K$ is a Dedekind domain, every integral ideal can be uniquely
expressed as a product of prime ideals.
We define the gcd of ideals $\a$ and $\b$ in the following way.\\
For
$$
\a = \prod_{\P\in \mathcal{P}_\K} \mathfrak{p}^{v_{ \mathfrak{p}}(\a)}, \phantom{m}
\b = \prod_{\P\in \mathcal{P}_\K} \mathfrak{p}^{v_{ \mathfrak{p}}(\b)}, 
$$
where all but finitely many $v_{\mathfrak{p}}(\a), v_{\mathfrak{p}}(\b)$ are zero, 
we define
$$
(\a,\b)=\gcd(\a,\b)=\prod_{\P\in \mathcal{P}_\K}\mathfrak{p}^{\min(v_{\mathfrak{p}}(\a),
v_{ \mathfrak{p}}(\b))},
$$
where we have denoted $\mathfrak{p}^{0}$ by $\mathcal{O}_\K$. 
Hence if prime factors of  $\a$ and $\b$ are all distinct,
$
(\a,\b) = \mathcal{O}_\K .
$
We notice that
$
(\a,\b)=\a+\b$
as
$v_\mathfrak{p}(\a+\b)=
\min (v_\mathfrak{p}(\a), v_\mathfrak{p}(\b)).$
From now on,
$ \mathcal{O}_\K$ will be denoted by 1.

\begin{lem}
For a number field $\K$ and a finite set $\Omega_f$, the limit
$$
\lim_{ x \to \infty} \left(  \frac{1}{\rho_\K} \sum_{ \substack{0\ne I \subset \O_{\K}
\atop{ \N(I) \leq x\atop  (I, P(\Omega_f))=1}}}
\frac{1}{\N(I)}  ~-~  \delta_\K(\Omega_f) \log x \right)
$$
exists and is denoted by $\gamma_\K(\Omega_f)$.
\end{lem}

\begin{proof}
Let $\Omega_f \subset \mathcal{P}_\K$ and $\P \in \mathcal{P}_\K$ be a prime ideal not in $\Omega_f$.
Using 
\begin{equation*}
\sum_{0\neq I\subseteq \O_\K \atop{\N(I)\leq x \atop{(I,\P P(\Omega_f) )=1}}}
\frac{1}{\N(I)}=
\sum_{0\neq I\subseteq \O_\K \atop{\N(I)\leq x \atop{(I,P(\Omega_f))=1}}}
\frac{1}{\N(I)}-\frac{1}{\N(\P)}
\sum_{0\neq I\subseteq \O_\K \atop{\N(I)\leq \frac{x}{\N(\P)} 
\atop{(I,P(\Omega_f))=1}}}\frac{1}{\N(I)},
\end{equation*}
the result follows by induction on cardinality of $\Omega_f$.
\end{proof}

\begin{lem}\label{finiteOmega}
Let $\Omega_f$ be a finite set of non-zero prime ideals. Then
$$
\gamma_\K(\Omega_f) = \delta_\K(\Omega_f)  \left(  \gamma_\K +  \sum_{ \mathfrak{p}\in 
\Omega_f} \frac {\log \N(\mathfrak{p})}{ \N(\mathfrak{p}) - 1} \right).
$$ 
\end{lem}

\begin{proof}
We have
\begin{align*}
\sum_{ \substack{0\ne I \subset \O_{\K}\atop{\N(I)\leq x\atop (I, P(\Omega_f)) = 1}}}
\frac{1}{\N(I)} &=  \sum_{ \substack{0\ne I \subset \O_{\K}\atop{ \N(I) \leq x}}}
\frac{1}{\N(I)} \sum_{J\mid(I,P(\Omega_f))} \mu(J)\\
& =\sum_{J\mid P(\Omega_f)} \frac{\mu(J)}{\N(J)} 
\sum_{ \substack{0\ne J_0 \subset \O_{\K}\atop{ \N(J_0) \leq \frac{x}{\N(J)}}}} 
\frac{1}{\N(J_0)}\\
&=\sum_{J \mid P(\Omega_f)} \frac{\mu(J)}{\N(J)}
\big\{\rho_\K \log \frac{x}{\N(J)} +\rho_{\K} \gamma_\K + \ro(1)\big\}\\
&= \delta_\K(\Omega_f) \big( \rho_\K \log x + \rho_{\K}\gamma_\K + \ro(1)\big)
 - \rho_\K\sum_{J\mid P(\Omega_f)} \frac{\mu(J)}{\N(J)} \log \N(J).
\end{align*}
We now consider the last term.
\begin{align*}
-\sum_{J\mid P(\Omega_f)} \frac{\mu(J)}{\N(J)} \log \N(J) 
&= \sum_{J \mid P(\Omega_f)} \frac{1}{\N(J)} \sum_{J_0\mid J} \Lambda(J_0)
\mu\left(JJ_0^{-1}\right)\\
&= \sum_{J_0 \mid P(\Omega_f)} \frac{\Lambda(J_0)}{ \N(J_0)}
\sum_{J_1\mid P(\Omega_f)J_0^{-1}} \frac{\mu(J_1)}{\N(J_1)}\\
& = \sum_{\mathfrak{p'} \in \Omega_f} 
\frac{\Lambda(\mathfrak{p'})}{ \N(\mathfrak{p'})}
\sum_{J_1\mid P(\Omega_f)\mathfrak{p'}^{-1}} 
\frac{\mu(J_1)}{\N(J_1)}\\
&= \sum_{\mathfrak{p'} \in \Omega_f} 
\frac{\log\N(\mathfrak{p'})}{ \N(\mathfrak{p'})}\left( 
\frac{\delta_\K(\Omega_f)}{1- \frac{1}{\N(\mathfrak{p'})}}\right)\\
& = \delta_\K(\Omega_f) \sum_{\mathfrak{p} \in \Omega_f} \frac {\log 
\N(\mathfrak{p})}{ \N(\mathfrak{p}) - 1}.
\end{align*}
Thus
\begin{align*}
\lim_{x\rightarrow \infty} \left(\frac{1}{\rho_\K} 
\sum_{ \substack{0\ne I \subset \O_{\K}\atop{\N(I)\leq x\atop (I, P(\Omega_f)) = 1}}}
\frac{1}{\N(I)} -  \delta_\K(\Omega_f)  \log x \right)
&= \delta_\K(\Omega_f)  \left(  \gamma_\K +  \sum_{ \mathfrak{p}\in \Omega_f} 
\frac {\log \N(\mathfrak{p})}{ \N(\mathfrak{p}) - 1} \right). 
\end{align*}
\end{proof}

\begin{cor}\label{infiniteOmega}
For a number field $\K$ and any set $\Omega \subset \mathcal{P}_\K$ satisfying \eqref{conv},
the limit
$$
\lim_{ x \to \infty} \left( \frac{1}{\rho_\K} 
\sum_{ \substack{0\ne I \subset \O_{\K}\atop{\N(I) \leq x\atop (I, P(\Omega(x)))=1}}}
\frac{1}{\N(I)}  ~-~  \delta_\K(\Omega(x)) \log x \right)
$$
exists and equals
$$ \delta_\K(\Omega)  \left(  \gamma_\K + \sum_{ \mathfrak{p}\in \Omega} 
\frac {\log \N(\mathfrak{p})}{ \N(\mathfrak{p}) - 1} \right).
$$ 
We denote this limit by $\gamma_\K(\Omega)$.
\end{cor}
\begin{proof}
Follows from \lemref{finiteOmega} since $\Omega(x)$ is a finite set.
\end{proof}

\section{Proof of \thmref{trans}}
Suppose there exist $i,j \in I$ such that
$$
\frac{\gamma_\K(\Omega_i)}{\delta_\K(\Omega_i)} \phantom{m} \text{and} \phantom{m}
\frac{\gamma_\K(\Omega_j)}{ \delta_\K(\Omega_j)}
$$
are algebraic.
Using \corref{infiniteOmega}, 
\begin{align}\label{star}
\frac{\gamma_\K(\Omega_i)}{\delta_\K(\Omega_i)} -  \frac{\gamma_\K(\Omega_j)}
{\delta_\K(\Omega_j)} &=  \sum_{ \mathfrak{p}\in \Omega_i} 
\frac {\log \N(\mathfrak{p})}{ \N(\mathfrak{p}) - 1} 
-\sum_{ \mathfrak{p}\in \Omega_j} 
\frac {\log \N(\mathfrak{p})}{ \N(\mathfrak{p}) - 1},
\end{align}
which is also an algebraic number. 
Since the sets $N_{\Omega_i} \backslash N_{\Omega_j}$ and 
$N_{\Omega_i} \backslash N_{\Omega_j}$ are non-empty and finite, 
the sets $\Omega_i \backslash  \Omega_j$ and $\Omega_j \backslash \Omega_i$ 
are also finite. Let
$$
\Omega_i \backslash \Omega_j= \{\mathfrak{p}_1, \mathfrak{p}_2 , 
\cdots, \mathfrak{p}_n\}, 
\phantom{mm} \Omega_j \backslash \Omega_i= \{\mathfrak{q}_1, \mathfrak{q}_2 , 
\cdots, \mathfrak{q}_m\}.
$$
Then equation \eqref{star} implies
\begin{equation}\label{UPF}
\sum_{ \mathfrak{p}\in \Omega_i} \frac {\log \N(\mathfrak{p})}{ \N(\mathfrak{p}) - 1} 
- \sum_{ \mathfrak{p}\in \Omega_j} \frac{\log \N(\mathfrak{p})}{\N(\mathfrak{p}) - 1}
= \sum_{s = 1}^{n} \frac {\log p^{f_s}_s}{ p^{f_s}_s - 1} - \sum_{t=1}^{m} 
\frac {\log q^{g_t}_t}{ q^{g_t}_t - 1}
= \log \left( \frac{\prod_{s=1}^{n} p^{\left(\frac{f_s}{p^{f_s}_s -1}\right)}_s}
{\prod_{t=1}^{m} q^{\left(\frac{g_t}{q^{g_t}_t -1}\right)}_t}\right),
\end{equation}
where  $\N(\mathfrak{p}_s)=p_s^{f_s}$ and $\N(\mathfrak{q}_t)=q_t^{g_t}$.
Using \lemref{FL} and unique prime factorisation of natural numbers, the expression
in \eqref{UPF} becomes a transcendental number, which gives a contradiction.

\bigskip
\noindent
{\bf Acknowledgments.} 
The authors would like to thank Prof. S. Gun for suggesting 
the problem and IMSc for providing academic facilities. They 
would also like to thank  Prof. P. Moree and J. Sivaraman for comments on an 
earlier version of the manuscript, which improved the 
exposition. The authors would like to thank the anonymous refree for 
careful reading of the paper and helpful comments.
The first author would like to thank the Max-
Planck-Institut f\"ur Mathematik for providing a friendly 
atmosphere.
The second author would like to thank Number Theory plan 
project, Department of Atomic Energy for financial support.


\smallskip
\begin{thebibliography}{xx}

\bibitem{DF}
H.  Diamond and K. Ford,
{\em Generalized Euler constants}, 
Math. Proc. Cambridge Philos. Soc. {\bf 145} (2008), no. 1, 27--41. 

\bibitem{LE}
L. Euler,
{\em De Progressionibus harmonicis observationes},
Commentarii academiae scientiarum Petropolitanae, {\bf 7} (1740), 150--161

\bibitem{YI2}
Y. Ihara,
{\em The Euler-Kronecker invariants in various families of global fields}, 
Arithmetics, geometry and coding theory (AGCT 2005), 79--102, 
Sémin. Congr., {\bf 21}, Soc. Math. France.

\bibitem{lag}
 Jeffrey C. Lagarias,
{\em Euler’s constant: Euler’s work and modern developments},
Bull. Amer. Math. Soc. {\bf 50} (2013), 527--628.

\bibitem{FL}
F. Lindemann,
{\em \"Uber die Zahl $\pi$},
Math. Ann. {\bf 20} (1882), no. 2, 213--225.

\bibitem{RM}
M. Ram Murty,
Problems in analytic number theory,
Second edition, {\em Springer}, New York, 2008.

\bibitem{ME}
M. Ram Murty and J. Esmonde,
Problems in Algebraic Number Theory,
{\em Springer-Verlag}, New York, 2005.

\bibitem{MZ}
M. Ram Murty and A. Zaytseva,
{\em Transcendence of generalized Euler constants},  
American Math. Monthly, {\bf 120} (2013), no. 1, 48--54.

\bibitem{MiRo}
M. Rosen,
{\em A generalization of Mertens' theorem},
J. Ramanujan Math. Soc. {\bf 14} (1999), no. 1, 1--19.

\end{thebibliography}
\end{document}